\input ./style/arxiv-general.cfg
\documentclass[aap,MSNbibl,seceqn,dvips]{arximspdf}
\makeatletter
   \@ifpackageloaded{graphicx}{}{\usepackage{graphicx}}
\makeatother

\doi{10.1214/15-AAP1128}
\volume{26}
\issue{3}
\pubyear{2016}
\firstpage{1636}
\lastpage{1658}
\docsubty{FLA}

\makeatletter
\newcommand{\mat}[4]{{\left(\matrix{ #1 & #2 \cr #3 & #4 }\right)}}
\newcommand{\eqref}[1]{(\ref{#1})}
\newcommand{\cov}{\operatorname{Cov}}
\newcommand{\Cor}{\operatorname{Cor}}
\newcommand{\E}{\mathbb{E}}
\newcommand{\R}{\mathbb{R}}
\newcommand{\N}{\mathbb{N}}
\newcommand{\p}{\mathbb{P}}
\newcommand{\B}{\mathcal{B}}
\newcommand{\D}{\mathcal{D}}
\newcommand{\Q}{\mathcal{Q}}
\newcommand{\C}{\mathcal{C}}
\newcommand{\V}{\mathcal{V}}
\newcommand{\T}{\mathcal{T}}
\renewcommand{\L}{\mathcal{L}}
\newcommand{\diag}{\operatorname{diag}}
\newcommand{\id}{\mathrm{I}}
\newcommand{\LN}{\operatorname{LN}}
\newcommand{\bd}{\mathbf}
\newtheorem{theorem}{Theorem}[section]
\newtheorem{corollary}[theorem]{Corollary}

\newtheorem{proposition}[theorem]{Proposition}
\newproclaim{definition}{Definition}[section]
\newproclaim{example}{Example}[section]

\newtheorem{question}{Question}
\newproclaim{remark}{Remark}[section]
\makeatother

\begin{document}
\begin{frontmatter}

\title{Bernoulli and tail-dependence compatibility}
\runtitle{Bernoulli and tail-dependence compatibility}

\begin{aug}
\author[A]{\fnms{Paul}~\snm{Embrechts}\thanksref{m1,T1}\ead[label=e1]{embrechts@math.ethz.ch}},
\author[B]{\fnms{Marius}~\snm{Hofert}\corref{}\thanksref{m2,T2}\ead[label=e2]{marius.hofert@uwaterloo.ca}}
\and
\author[B]{\fnms{Ruodu}~\snm{Wang}\thanksref{m2,T2,T3}\ead[label=e3]{wang@uwaterloo.ca}}
\thankstext{T1}{Supported by the Swiss Finance Institute.}
\thankstext{T2}{Supported by the Natural Sciences and Engineering Research Council of Canada
(NSERC Discovery Grant numbers 5010 and 435844, resp.).}
\thankstext{T3}{Supported
by the Forschungsinstitut f\"ur Mathematik (FIM) at ETH Zurich.}
\runauthor{P. Embrechts, M. Hofert and R. Wang}
\affiliation{ETH Zurich\thanksmark{m1} and University of
Waterloo\thanksmark{m2}}
\address[A]{P. Embrechts\\
RiskLab\\
Department of Mathematics\\
\quad and Swiss Finance Institute\\
ETH Zurich\\
8092 Zurich\\
Switzerland\\
\printead{e1}}

\address[B]{M. Hofert\\
R. Wang\\
Department of Statistics\\
\quad and Actuarial Science\\
University of Waterloo\\
200 University Avenue West\\
Waterloo, Ontario N2L 3G1\\
Canada\\
\printead{e2}\\
\phantom{E-mail:\ }\printead*{e3}}
\end{aug}

%
\received{\smonth{1} \syear{2015}}

%
\begin{abstract}
The tail-dependence compatibility problem is introduced. It raises the
question whether a given $d\times d$-matrix of entries in the unit
interval is
the matrix of pairwise tail-dependence coefficients of a $d$-dimensional
random vector. The problem is studied together with Bernoulli-compatible
matrices, that is, matrices which are expectations of outer products of random
vectors with Bernoulli margins. We show that a square matrix with diagonal
entries being 1 is a tail-dependence matrix if and only if it is a
Bernoulli-compatible matrix multiplied by a constant. We introduce new copula
models to construct tail-dependence matrices, including commonly used matrices
in statistics.
\end{abstract}

%
\begin{keyword}[class=AMS]
\kwd{60E05}
\kwd{62H99}
\kwd{62H20}
\kwd{62E15}
\kwd{62H86}
\end{keyword}
\begin{keyword}
\kwd{Tail dependence}
\kwd{Bernoulli random vectors}
\kwd{compatibility}
\kwd{matrices}
\kwd{copulas}
\kwd{insurance application}
\end{keyword}
%
\end{frontmatter}

\section{Introduction}
The problem of how to construct a bivariate random vector $(X_1,X_2)$ with
log-normal marginals $X_1\sim\LN(0,1)$, $X_2\sim\LN(0,16)$ and
correlation coefficient
$\Cor(X_1,X_2)=0.5$ is well known in the history of dependence modeling,
partially because of its relevance to risk management practice. The
short answer
is: There is no such model; see Embrechts et al. \cite{EMS02} who
studied these kinds of
problems in terms of copulas. Problems of this kind were brought to
RiskLab at
ETH Zurich by the insurance industry in the mid-1990s when dependence was
thought of in terms of correlation (matrices). For further background on
quantitative risk management, see McNeil et al. \cite{MFE15}. Now,
almost 20 years later,
copulas are a well established tool to quantify dependence in
multivariate data
and to construct new multivariate distributions. Their use has become standard
within industry and regulation. Nevertheless, dependence is still
summarized in
terms of numbers [as opposed to (copula) functions], so-called \emph
{measures of
association}. Although there are various ways to compute such numbers in
dimension $d>2$, measures of association are still most widely used in the
bivariate case $d=2$. A~popular measure of association is tail
dependence. It is
important for applications in quantitative risk management as it
measures the
strength of dependence in either the lower-left or upper-right tail of the
bivariate distribution, the regions quantitative risk management is mainly
concerned with.

We were recently asked\setcounter{footnote}{3}\footnote{By Federico Degen (Head Risk Modeling and
Quantification, Zurich Insurance Group) and Janusz Milek (Zurich
Insurance Group).} the following question which is in
the same spirit as the log-normal correlation problem if one replaces
``correlation'' by ``tail dependence''; see Section~\ref{sec:Lambda:1}
for a definition.

\begin{quote}
\textit{For which $\alpha\in[0,1]$ is the matrix}
%
\begin{equation}
\Gamma_d(\alpha)=\pmatrix{ 1 & 0 & \cdots& 0 & \alpha\vspace*{2pt}
\cr
0 & 1 & \cdots& 0 & \alpha\vspace*{2pt}
\cr
\vdots& \vdots& \ddots& \vdots&
\vdots\vspace*{2pt}
\cr
0 & 0 & \cdots& 1 & \alpha\vspace*{2pt}
\cr
\alpha& \alpha&
\cdots& \alpha& 1 }\label{eq:Gammad}
\end{equation}
\textit{a matrix of pairwise (either lower or upper) tail-dependence coefficients?}
\end{quote}

Intrigued by this question, we more generally consider the following
\emph{tail-dependence compatibility problem} in this paper:

\begin{quote}
\textit{When is a given matrix in $[0,1]^{d\times d}$ the
matrix of pairwise (either lower or upper) tail-dependence
coefficients?}
\end{quote}

In what follows, we call a matrix of pairwise tail-dependence
coefficients a
\emph{tail-dependence matrix}. The compatibility problems of tail-dependence
coefficients were studied in \cite{J97}. In particular, when $d=3$,
inequalities
for the bivariate tail-dependence coefficients have been established; see
Joe \cite{J97}, Theorem~3.14, as well as Joe \cite{J14}, Theorem~8.20.
The sharpness of these inequalities is obtained in
\cite{NJL09}. It is generally open to characterize the tail-dependence matrix
compatibility for $d>3$.

Our aim in this paper is to give a full answer to the tail-dependence
compatibility problem; see Section~\ref{sec:Lambda}. To this end, we introduce
and study \emph{Bernoulli-compatible matrices} in Section~\ref
{sec:bern}. As a
main result, we show that a matrix with diagonal entries being 1 is a compatible
tail-dependence matrix if and only if it is a Bernoulli-compatible matrix
multiplied by a constant. In Section~\ref{sec:model}, we provide probabilistic
models for a large class of tail-dependence matrices, including
commonly used
matrices in statistics. Section~\ref{sec:con} concludes.

Throughout this paper, $d$ and $m$ are positive integers, and we
consider an
atomless probability space $(\Omega, \mathcal A, \p)$ on which all random
variables and random vectors are defined. Vectors are considered as column
vectors. For two matrices $A,B$, $B\ge A$ and $B\le A$ are understood as
component-wise inequalities. We let $A\circ B$ denote the Hadamard product,
that is, the element-wise product of two matrices $A$ and $B$ of the same
dimension. The $d\times d$ identity matrix is denoted by $I_d$. For a
square matrix
$A$, $\diag(A)$ represents a diagonal matrix with diagonal entries
equal to
those of $A$, and $A^\top$ is the transpose of $A$. We denote $\id_E$ the
indicator function of an event (random or deterministic) $E\in\mathcal A$.
$\mathbf{0}$ and $\mathbf{1}$ are vectors with all components being 0
and 1, respectively, as long as the dimension of the vectors is clear
from the context.

\section{Bernoulli compatibility}\label{sec:bern}
In this section, we introduce and study the \emph{Bernoulli-compatibility
problem}. The results obtained in this section are the basis for the
\emph{tail-dependence compatibility problem} treated in
Section~\ref{sec:Lambda}; many of them are of independent interest,
for example, for the
simulation of sequences of Bernoulli random variables.

\subsection{Bernoulli-compatible matrices}
%
\begin{definition}[(Bernoulli vector, $\V_d$)]
A \emph{Bernoulli vector} is a random vector $\mathbf{X}$ supported
by $\{0,1\}^d$ for
some $d\in\N$. The set of all $d$-Bernoulli vectors is denoted by $\V_d$.
\end{definition}
Equivalently, $\mathbf{X}=(X_1,\ldots,X_d)$ is a Bernoulli vector if
and only if
$X_i\sim\mathrm{B}(1,p_i)$ for some $p_i\in[0,1]$, $i=1,\ldots,d$.
Note that
here we do not make any assumption about the dependence structure among the
components of $\mathbf{X}$. Bernoulli vectors play an important role
in credit risk
analysis; see, for example, Bluhm and Overbeck \cite{BO06} and Bluhm
et al. \cite{BOW02}, Section~2.1.

In this section, we investigate the following question which we refer
to as
the \emph{Bernoulli-compatibility problem}.
%
\begin{question}\label{Q:B:compatibility}
Given a matrix $B\in[0,1]^{d\times d}$, can we find a Bernoulli vector
$\mathbf{X}$
such that $B=\E[\mathbf{X}\mathbf{X}^\top]$?
\end{question}
For studying the Bernoulli-compatibility problem, we introduce the
notion of
Bernoulli-compatible matrices.
%
\begin{definition}[(Bernoulli-compatible matrix, $\B_d$)]
A $d\times d$ matrix $B$ is a \emph{Bernoulli-compatible matrix}, if
$B=\E[\mathbf{X}\mathbf{X}^\top]$ for some $\mathbf{X}\in\V_d$.
The set of all $d\times d$
Bernoulli-compatible matrices is denoted by $\B_d$.
\end{definition}
Concerning covariance matrices, there is extensive research on the compatibility
of covariance matrices of Bernoulli vectors in the realm of statistical
simulation and time series analysis; see, for example, Chaganty and Joe
\cite{CJ06}. It is known that,
when $d\ge3$, the set of all compatible $d$-Bernoulli correlation
matrices is
strictly contained in the set of all correlation matrices. Note that
$\E[\mathbf{X}\mathbf{X}^\top]=\cov(\mathbf{X})+\E[\mathbf
{X}]\E[\mathbf{X}]^\top$. Hence, Question~\ref{Q:B:compatibility}
is closely
related to the characterization of compatible Bernoulli covariance matrices.

Before we characterize the set $\B_d$ in Section~\ref{chara:Bcm}, and thus
address Question~\ref{Q:B:compatibility}, we first collect some facts
about elements of $\B_d$.
%
\begin{proposition}\label{prop:prelim}
Let $B,B_1,B_2\in\B_d$. Then:
\begin{longlist}[(iii)]
\item[(i)]$B\in[0,1]^{d\times d}$.
\item[(ii)] $\max\{b_{ii}+b_{jj}-1,0\}\le b_{ij}\le
\min\{b_{ii},b_{jj}\}$ for $i,j=1,\ldots,d$ and $B=(b_{ij})_{d\times d}$.
\item[(iii)]$t B_1+(1-t)B_2 \in\B_d$ for $t\in[0,1]$, that is, $\B_d$ is
a convex
set.
\item[(iv)] $B_1\circ B_2 \in\B_d$, that is, $\B_d$
is closed under the Hadamard product.
\item[(v)]$(0)_{d\times d}\in\B_d$ and $(1)_{d\times d}\in\B_d$.
\item[(vi)] For any $\bd p=(p_1,\ldots,p_d)\in[0,1]^d$, the matrix
$B=(b_{ij})_{d\times d}\in\B_d$ where $b_{ij}=p_i p_j$ for $i\ne j$
and $b_{ii}=p_i$, $i,j=1,\ldots,d$.
\end{longlist}
\end{proposition}
\begin{pf}
Write $B_1=\E[\mathbf{X}\mathbf{X}^\top]$ and $B_2=\E[\mathbf
{Y}\mathbf{Y}^\top]$ for
$\mathbf{X},\mathbf{Y}\in\V_d$, and $\mathbf{X}$ and $\mathbf{Y}$
are independent.
\begin{longlist}[(iii)]
\item[(i)] Clear.
\item[(ii)] This directly follows from the Fr\'echet--Hoeffding bounds; see
McNeil et al. \cite{MFE15}, Remark~7.9. 
\item[(iii)] Let $A\sim\mathrm{B}(1,t)$ be a
Bernoulli random variable independent of $\mathbf{X},\mathbf{Y}$, and let
$\mathbf{Z}=A\mathbf{X}+(1-A)\mathbf{Y}$. Then $\mathbf{Z}\in\V
_d$, and $\E[\mathbf{Z}\mathbf{Z}^\top]=t
\E[\mathbf{X}\mathbf{X}^\top]+(1-t)\E[\mathbf{Y}\mathbf{Y}^\top
]=t B_1+(1-t)B_2$. Hence, $t B_1+(1-t)B_2\in
\B_d$.
\item[(iv)] Let $\bd p=(p_1,\ldots,p_d)$, $\bd q=(q_1,\ldots,q_d)\in\R^d$. Then
\begin{eqnarray*}
(\bd p\circ\bd q) (\bd p\circ\bd q)^\top&=&(p_iq_i)_{d}
(p_iq_i)_d^\top
=(p_iq_ip_jq_j)_{d\times d}=(p_ip_j)_{d\times d}
\circ(q_iq_j)_{d\times
d}\\
&=&\bigl(\bd p\bd
p^\top\bigr)\circ\bigl(\bd q\bd q^\top\bigr).
\end{eqnarray*}
Let $\mathbf{Z}=\mathbf{X}\circ\mathbf{Y}$. It follows that
$\mathbf{Z}\in\V_d$ and $\E[\mathbf{Z}\mathbf{Z}^\top]=\E
[(\mathbf{X}\circ\mathbf{Y}) (\mathbf{X}\circ\mathbf{Y})^\top
]=\E[(\mathbf{X}\mathbf{X}^\top)\circ(\mathbf{Y}\mathbf{Y}^\top
)]=\E[\mathbf{X}\mathbf{X}^\top]\circ\E[\mathbf{Y}\mathbf
{Y}^\top]=B_1\circ B_2$. Hence, $B_1\circ B_2 \in\B_d$.
\item[(v)] Consider $\mathbf{X}=\mathbf{0}\in\V_d$. Then $(0)_{d\times
d}=\E[\mathbf{X}\mathbf{X}^\top]\in\B_d$ and similarly for
$(1)_{d\times d}$.
\item[(vi)] Consider $\mathbf{X}\in\V_d$ with independent components and
$\E[\mathbf{X}]=\bd p$.\quad\qed
\end{longlist}
\noqed\end{pf}

\subsection{Characterization of Bernoulli-compatible matrices}\label
{chara:Bcm}
We are now able to give a characterization of the set $\B_d$ of
Bernoulli-compatible matrices and thus address Question~\ref
{Q:B:compatibility}.
%
\begin{theorem}[(Characterization of $\B_d$)]\label{thm:bern}
$\B_d$ has the following characterization:
%
\begin{equation}\label{bern-1}
\B_d= \Biggl\{ \sum_{i=1}^{n}
a_i \bd{p}_i\bd{p}_i^\top:
\bd{p}_i\in\{0,1\}^d, a_i\ge0, i=1,\ldots,n,
\sum_{i=1}^n a_i=1, n\in \N
\Biggr\};\hspace*{-25pt}
\end{equation}
that is, $\B_d$ is the convex hull of $\{\bd p\bd p^\top: \bd p\in\{
0,1\}^d\}$.
In particular, $\B_d$ is closed under convergence in the Euclidean norm.
\end{theorem}
\begin{pf}
Denote the right-hand side of \eqref{bern-1} by $\mathcal M$. For
$B\in
\B_d$, write $B=\E[\mathbf{X}\mathbf{X}^\top]$ for some $\mathbf
{X}\in\V_d$. It follows that
\[
B=\sum_{\bd p \in\{0,1\}^d}\bd p \bd p^\top\p(\mathbf{X}=
\bd p)\in \mathcal M,
\]
hence $\B_d \subseteq\mathcal M$. Let $\mathbf{X}=\mathbf{p}\in\{
0,1\}^d$. Then $\mathbf{X}\in\V_d$ and $\E[\mathbf{X}\mathbf
{X}^\top]=\bd p
\bd p^\top\in\B_d$. By Proposition~\ref{prop:prelim}, $\B_d$ is a
convex set
which contains $\{\bd p\bd p^\top: \bd p\in\{0,1\}^d\}$, hence
$\mathcal
M\subseteq\B_d$. In summary, $\mathcal M=\B_d$. From \eqref
{bern-1}, we can
see that $\B_d$ is closed under convergence in the Euclidean norm.
\end{pf}

A matrix $B$ is \emph{completely positive} if $B=AA^\top$ for some (not
necessarily square) matrix $A\ge0$. Denote by $\C_d$ the set of completely
positive matrices. It is known that $\C_d$ is the convex cone with extreme
directions $\{\bd p\bd p^\top: \bd p\in[0,1]^d\}$; see, for example,
R\"uschendorf \cite{R81} and Berman and Shaked-Monderer \cite{BS03}.
We thus obtain the following result.
%
\begin{corollary}\label{cpm}
Any Bernoulli-compatible matrix is completely positive.
\end{corollary}

\begin{remark}\label{rem:Bd}
One may wonder whether $B=\E[\mathbf{X}\mathbf{X}^\top]$ is
sufficient to determine the
distribution of $\mathbf{X}$, that is, whether the
decomposition
%
\begin{equation}
B=\sum_{i=1}^{2^d} a_i
\bd{p}_i\bd{p}_i^\top\label{op1}
\end{equation}
is unique for distinct vectors $\bd{p}_i$ in $\{0,1\}^d$. While the
decomposition is trivially unique for $d=2$, this is in general
false for $d\ge3$, since there are $2^d-1$ parameters in \eqref{op1}
and only
$d(d+1)/2$ parameters in $B$. The following is an example for $d=3$.
Let
\begin{eqnarray*}
B&=&\frac{1}4\pmatrix{ 2 & 1 & 1 \vspace*{2pt}
\cr
1 & 2 & 1 \vspace*{2pt}
\cr
1 & 1 & 2 }
\\
&=&\frac{1}4 \bigl( (1,1,1)^\top(1,1,1) + (1,0,0)^\top(1,0,0)
+ (0,1,0)^\top(0,1,0)\\
&&{} + (0,0,1)^\top(0,0,1) \bigr)
\\
&=&\frac{1}4 \bigl( (1,1,0)^\top(1,1,0) + (1,0,1)^\top(1,0,1)
+ (0,1,1)^\top(0,1,1)\\
&&{} + (0,0,0)^\top(0,0,0) \bigr).
\end{eqnarray*}
Thus, by combining the above two decompositions, $B\in\B_3$ has infinitely
many different decompositions of the form \eqref{op1}. Note that, as
in the case of
completely positive matrices, it is generally difficult to find
decompositions of
form \eqref{op1} for a given matrix $B$.
\end{remark}




\subsection{Convex cone generated by Bernoulli-compatible matrices}
In this section, we study the convex cone generated by $\B_d$, denoted
by $\B_d^*$:
%
\begin{equation}
\B_d^*=\{aB: a\ge0, B\in\B_d\}.\label{bern-2}
\end{equation}
The following proposition is implied by Proposition~\ref{prop:prelim}
and Theorem~\ref{thm:bern}.
%
\begin{proposition}\label{prop:bern}
$\B_d^*$ is the convex cone with extreme directions $\{\bd p\bd p^\top
: \bd
p\in\{0,1\}^d\}$. Moreover, $\B_d^*$ is a commutative semiring
equipped with
addition $(\B_d^*, +)$ and multiplication $(\B_d^*, \circ)$.
\end{proposition}
It is obvious that $\B_d^*\subseteq\C_d$. One may wonder whether $\B
_d^*$ is
identical to $\C_d$, the set of completely positive matrices. As the
following example shows, this is false in general for $d\ge2$.
%
\begin{example}
Note that $B \in\B_d^*$ also satisfies Proposition~\ref{prop:prelim},
part~(ii). Now consider $\bd p=(p_1,\ldots,p_d)\in
(0,1)^d$ with
$p_i> p_j$ for some $i\neq j$. Clearly, $\bd p\bd p^\top\in\C_d$,
but $p_ip_j>
p_j^2=\min\{p_i^2,p_j^2\}$ contradicts Proposition~\ref{prop:prelim},
part~(ii), hence
$\bd p\bd p^\top\notin\B_d^*$.
\end{example}

For the following result, we need the notion of diagonally dominant
matrices. A
matrix $A\in\R^{d\times d}$ is called \emph{diagonally dominant} if,
for all
$i=1,\ldots,d$, $\sum_{j\neq i}|a_{ij}|\le|a_{ii}|$.
%
\begin{proposition}\label{ddm}
Let $\D_d$ be the set of nonnegative, diagonally dominant $d\times
d$-matrices. Then $\D_d\subseteq\B_d^*$.
\end{proposition}
\begin{pf}
For $i,j=1,\ldots,d$, let $\bd p^{(ij)}=(p^{(ij)}_1,\ldots,p^{(ij)}_d)$
where $p^{(ij)}_k=\id_{\{k=i\}\cup\{k=j\}}$. It is straightforward to
verify that the
$(i,i)$-, $(i,j)$-, $(j,i)$- and $(j,j)$-entries of the matrix
$M^{(ij)}=\bd
p^{(ij)} (\bd p^{(ij)})^\top$ are 1, and the other entries are 0. For
$D=(d_{ij})_{d\times d}\in\D_d$, let
\[
D^*=\bigl(d_{ij}^*\bigr)_{d\times d}=\sum
_{i=1}^d\sum_{j=1, j\ne i}^d
d_{ij}M^{(ij)}.
\]
By Proposition~\ref{prop:bern}, $D^*\in\B_d^*$. It follows that
$d_{ij}^*=d_{ij}$ for $i\ne j$ and $d_{ii}^*=\sum_{j=1,j\ne i}^d
d_{ij}\le
d_{ii}$. Therefore, $D=D^*+\sum_{i=1}^d (d_{ii}-d_{ii}^*)M^{(ii)}$, which,
by Proposition~\ref{prop:bern}, is in $\B_d^*$.
\end{pf}

For studying the tail-dependence compatibility problem in Section~\ref
{sec:Lambda}, the subset
\[
\B^I_d=\bigl\{B: B\in\B^*_d,
\diag(B)=I_d\bigr\}
\]
of $\B_d^*$ is of interest. It is straightforward to see from
Proposition~\ref{prop:prelim}
and Theorem~\ref{thm:bern} that $\B^I_d$ is a convex set, closed
under the Hadamard product and
convergence in the Euclidean norm. These properties of $\B^I_d$ will
be used later.

\section{Tail-dependence compatibility}\label{sec:Lambda}
\subsection{Tail-dependence matrices}\label{sec:Lambda:1}
The notion of tail dependence captures (extreme) dependence in the
lower-left or
upper-right tails of a bivariate distribution. In what follows, we
focus on
lower-left tails; the problem for upper-right tails follows by a reflection
around $(1/2,1/2)$, that is, studying the survival copula of the
underlying copula.
%
\begin{definition}[(Tail-dependence coefficient)]
The \emph{(lower) tail-dependence coefficient} of two continuous
random variables $X_1\sim F_1$ and $X_2\sim F_2$
is defined by
%
\begin{equation}
\lambda=\lim_{u\downarrow0}\frac{\p(F_1(X_1)\le u, F_2(X_2)\le
u)}{u},\label{def:tail:dep}
\end{equation}
given that the limit exists.
\end{definition}
If we denote the copula of $(X_1,X_2)$ by $C$, then
\[
\lambda=\lim_{u\downarrow0}\frac{C(u,u)}{u}.
\]
Clearly, $\lambda\in[0,1]$, and $\lambda$ only depends on the copula of
$(X_1,X_2)$, not the marginal distributions. For virtually all copula models
used in practice, the limit in \eqref{def:tail:dep} exists; for how to
construct
an example where $\lambda$ does not exist; see Kortschak and
Albrecher \cite{kortschakalbrecher2009}.
%
\begin{definition}[(Tail-dependence matrix, $\T_d$)]
Let $\mathbf{X}=(X_1,\ldots,X_d)$ be a random vector with continuous marginal
distributions. The \emph{tail-dependence matrix} of $\mathbf{X}$ is
$\Lambda=(\lambda_{ij})_{d\times d}$, where $\lambda_{ij}$ is the
tail-dependence coefficient of $X_i$ and $X_j$, $i,j=1,\ldots,d$. We
denote by
$\T_d$ the set of all tail-dependence matrices.
\end{definition}
The following proposition summarizes basic properties of tail-dependence
matrices. Its proof is very similar to that of Proposition~\ref
{prop:prelim} and
is omitted here.
%
\begin{proposition}
For any $\Lambda_1, \Lambda_2\in\T_d$, we have that:
\begin{longlist}[(iii)]
\item[(i)]$\Lambda_1=\Lambda_1^\top$.
\item[(ii)]$t \Lambda_1+(1-t)\Lambda_2 \in\T_d$ for $t\in[0,1]$, that
is, $\T_d$ is
a convex set.
\item[(iii)]$I_d\le\Lambda_1\le(1)_{d\times d}$ with $I_d\in\T_d$ and
$(1)_{d\times d}\in\T_d$.
\end{longlist}
\end{proposition}
As we will show next, $\T_d$ is also closed under the Hadamard product.
%
\begin{proposition}\label{prop:hadamard:prod}
Let $k\in\N$ and $\Lambda_1,\ldots,\Lambda_k\in\T_d$. Then
$\Lambda_1\circ\cdots\circ\Lambda_k\in\T_d$.
\end{proposition}
\begin{pf}
Note that it would be sufficient to show the result for $k=2$, but we provide
a general construction for any $k$. For each $l=1,\ldots,k$, let $C_l$
be a
$d$-dimensional copula with tail-dependence matrix $\Lambda_l$. Furthermore,
let $g(u)=u^{1/k}$, $u\in[0,1]$. It follows from Liebscher \cite{L08} that
$C(u_1,\ldots,u_d)=\prod_{l=1}^k C_l(g(u_1),\ldots,g(u_d))$ is a
copula; note
that
%
\begin{equation}
\Bigl(g^{-1}\Bigl(\max_{1\le l\le
k}\{U_{l1}\}
\Bigr),\ldots,g^{-1}\Bigl(\max_{1\le l\le
k}\{U_{ld}
\}\Bigr) \Bigr)\sim C\label{eq:liebscher:cop}
\end{equation}
for independent random vectors $(U_{l1},\ldots,U_{ld})\sim C_l$,
$l=1,\ldots,k$. The $(i,j)$-entry $\lambda_{ij}$ of $\Lambda$
corresponding to $C$ is thus given by
\begin{eqnarray*}
\lambda_{ij}&=&\lim_{u\downarrow0}\frac{\prod_{l=1}^k
C_{l,ij}(g(u),g(u))}{u}=\lim
_{u\downarrow
0}\prod_{l=1}^k
\frac{C_{l,ij}(g(u),g(u))}{g(u)}\\
&=&\prod_{l=1}^k\lim
_{u\downarrow
0}\frac{C_{l,ij}(g(u),g(u))}{g(u)}
\\
&=&\prod_{l=1}^k\lim_{u\downarrow0}
\frac{C_{l,ij}(u,u)}{u}=\prod_{l=1}^k
\lambda_{l,ij},
\end{eqnarray*}
where $C_{l,ij}$ denotes the $(i,j)$-margin of $C_l$ and $\lambda_{l,ij}$
denotes the $(i,j)$th entry of $\Lambda_l$, $l=1,\ldots,k$.
\end{pf}

\subsection{Characterization of tail-dependence matrices}
In this section, we investigate the following question.
%
\begin{question}\label{Q:characterization}
Given a $d\times d$ matrix $\Lambda\in[0,1]^{d\times d}$, is
it a tail-dependence matrix? 
\end{question}
The following theorem fully characterizes tail-dependence matrices, and thus
provides a theoretical (but not necessarily practical) answer to
Question~\ref{Q:characterization}.
%
\begin{theorem}[(Characterization of $\T_d$)]\label{thm:main:characterization}
A square matrix with diagonal entries being 1 is a tail-dependence
matrix if
and only if it is a Bernoulli-compatible matrix multiplied by a
constant. Equivalently, $\T_d=\B^I_d$.
\end{theorem}
\begin{pf}
We first show that $\T_d\subseteq\B^I_d$. For each
$\Lambda=(\lambda_{ij})_{d\times d}\in\T_d$, suppose that $C$ is a copula
with tail-dependence matrix $\Lambda$ and $\mathbf{U}=(U_1,\ldots
,U_n)\sim C$. Let $\mathbf{W}_u=(\id_{\{U_1\le
u\}},\ldots,\id_{\{U_d\le u\}})$. By definition,
\[
\lambda_{ij}=\lim_{u\downarrow0}\frac{1}u\E[
\id_{\{U_i\le u\}}\id _{\{U_j\le u\}}]
\]
and
\[
\Lambda=\lim_{u\downarrow0}\frac{1}u \E\bigl[
\mathbf{W}_u\mathbf {W}_u^\top\bigr].
\]
Since $\B^I_d$ is closed and $\E[\mathbf{W}_u\mathbf{W}_u^\top
]/u\in\B^I_d$, we have that $\Lambda\in\B^I_d$.

Now consider $\B^I_d\subseteq\T_d$. By definition of $\B^I_d$, each
$B\in
\B^I_d$ can be written as $B=\E[\mathbf{X}\mathbf{X}^\top]/p$ for
an $\mathbf{X}\in\V_d$ and
$\E[\mathbf{X}]=(p,\ldots,p)\in(0,1]^d$. Let $U,V\sim\mathrm
{U}[0,1]$, $U,V,\mathbf{X}$ be
independent and
%
\begin{equation}
\mathbf{Y}=\mathbf{X}pU+(\mathbf{1}-\mathbf{X}) \bigl(p+(1-p)V
\bigr).\label{smodel}
\end{equation}
We can verify that for $t\in[0,1]$ and $i=1,\ldots,d$,
\begin{eqnarray*}
\p(Y_i\le t)&=&\p(X_i=1)\p(pU \le t)+
\p(X_i=0)\p\bigl(p+(1-p)V\le t\bigr)
\\
&=&p\min\{t/p,1\} +(1-p) \max\bigl\{(t-p)/(1-p),0\bigr\}=t,
\end{eqnarray*}
that is, $Y_1,\ldots,Y_d$ are $\mathrm{U}[0,1]$-distributed. Let
$\lambda_{ij}$
be the tail-dependence coefficient of $Y_i$ and $Y_j$, $i,j=1,\ldots,d$.
For $i,j=1,\ldots,d$ we obtain that
\begin{eqnarray*}
\lambda_{ij}&=&\lim_{u\downarrow0}\frac{1}u
\p(Y_i\le u, Y_j\le u)=\lim_{u\downarrow0}
\frac{1}u\p(X_i=1, X_j=1)\p(pU\le u)\\
&=&
\frac{1}p \E[X_iX_j].
\end{eqnarray*}
As a consequence, the tail-dependence matrix of $(Y_1,\ldots,Y_d)$ is
$B$ and $B\in\T_d$.
\end{pf}
It follows from Theorem~\ref{thm:main:characterization} and
Proposition~\ref{prop:bern} that $\T_d$ is
the ``1-diagonals'' cross-section of the convex cone with extreme directions
$\{\bd p\bd p^\top: \bd p\in\{0,1\}^d\}$. Furthermore, the proof of
Theorem~\ref{thm:main:characterization} is constructive. As we saw,
for any $B\in\B^I_d$, $\mathbf{Y}$
defined by \eqref{smodel} has tail-dependence matrix $B$. This interesting
construction will be applied in Section~\ref{sec:model} where we show that
commonly applied matrices in statistics are tail-dependence
matrices and where we derive the copula of $\mathbf{Y}$.
%
\begin{remark}
From the fact that $\T_d=\B_d^I$ and $\B_d^I$ is closed under the Hadamard
product [see Proposition~\ref{prop:prelim}, part~(iv)],
Proposition~\ref{prop:hadamard:prod} directly follows. Note, however,
that our
proof of Proposition~\ref{prop:hadamard:prod} is constructive. Given
tail-dependence matrices and corresponding copulas, we can construct a copula
$C$ which has the Hadamard product of the tail-dependence matrices as
corresponding tail-dependence matrix. If sampling of all involved
copulas is
feasible, we can sample $C$; see Figure~\ref{fig:liebscher} for
examples.\footnote{All
plots can be reproduced via the \textsf{R} package \texttt{copula} (version
$\ge$ 0.999-13) by calling \texttt{demo(tail\_compatibility)}.}

\begin{figure}

\includegraphics{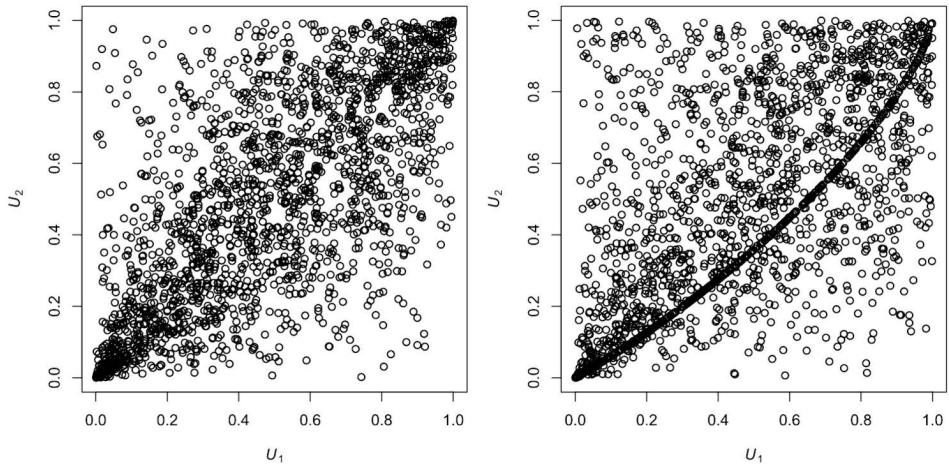}

\caption{Left-hand side: Scatter plot of 2000 samples from
\protect\eqref{eq:liebscher:cop} for $C_1$ being a Clayton copula
with parameter
$\theta=4$ ($\lambda_1=2^{-1/4}\approx0.8409$) and $C_2$ being a $t_3$
copula with parameter $\rho=0.8$ [tail-dependence coefficient
$\lambda_2=2t_4(-2/3)\approx0.5415$]. By Proposition~\protect\ref
{prop:hadamard:prod},
the tail-dependence coefficient of \protect\eqref{eq:liebscher:cop}
is thus
$\lambda=\lambda_1\lambda_2=2^{3/4}t_4(-2/3)\approx0.4553$. Right-hand
side: $C_1$ as before, but $C_2$ is a survival Marshall--Olkin copula with
parameters $\alpha_1=2^{-3/4},\alpha_2=0.8$, so that
$\lambda=\lambda_1\lambda_2=1/2$.}
\label{fig:liebscher}
\end{figure}
\end{remark}

Theorem~\ref{thm:main:characterization} combined with Corollary~\ref{cpm}
directly leads to the following result.
%
\begin{corollary}
Every tail-dependence matrix is completely positive, and hence positive
semi-definite.
\end{corollary}
Furthermore, Theorem~\ref{thm:main:characterization} and
Proposition~\ref{ddm} imply the
following result.

\begin{corollary}\label{coro:ddm}
Every diagonally dominant matrix with nonnegative entries and diagonal
entries being $1$ is a tail-dependence matrix.
\end{corollary}

Note that this result already yields the if-part of Proposition~\ref
{prop:example} below.

\section{Compatible models for tail-dependence matrices} \label{sec:model}
\subsection{Widely known matrices}
We now consider the following three types of matrices
$\Lambda=(\lambda_{ij})_{d\times d}$ which are frequently applied in
multivariate statistics and time series analysis and show that they are
tail-dependence matrices.
\begin{longlist}[(a)]
\item[(a)]
Equicorrelation matrix with parameter $\alpha
\in[0,1]$:
$\lambda_{ij}=\id_{\{i=j\}} +\alpha\id_{\{i\ne j\}}$, $ i,j=1,\ldots,d$.
\item[(b)]
AR(1) matrix with parameter $\alpha\in[0,1]$:
$\lambda_{ij}= \alpha^{ |i-j| }$,
$ i,j=1,\ldots,d$.
\item[(c)]
MA(1) matrix with parameter $\alpha\in
[0,1/2]$: $\lambda_{ij}=\id_{\{i=j\}}
+\alpha\id_{\{|i-j|=1\}}$, $ i,j=1,\ldots,d$.
\end{longlist}
Chaganty and Joe \cite{CJ06} considered the compatibility of
correlation matrices of
Bernoulli vectors for the above three types of matrices and obtained
necessary and sufficient conditions for the existence of compatible
models for $d=3$. For the tail-dependence compatibility problem that we
consider in this paper, the above three types of matrices are all
compatible, and we are able to construct corresponding models for each case.

\begin{proposition}
Let $\Lambda$ be the tail-dependence matrix of the $d$-dimension\-al
random vector
%
\begin{equation}
\mathbf{Y}=\mathbf{X}pU+(\mathbf{1}-\mathbf{X}) \bigl(p+(1-p)V
\bigr),\label{smodel2}
\end{equation}
where $U,V\sim\mathrm{U}[0,1]$, $\mathbf{X}\in\V_d$ and
$U,V,\mathbf{X}$ are independent.
\begin{longlist}[(iii)]
\item[(i)] For $\alpha\in[0,1]$, if $\mathbf{X}$ has independent
components and
$\E[X_1]=\cdots=\E[X_d]=\alpha$, then $\Lambda$ is an
equicorrelation matrix with
parameter $\alpha$; that is, \textup{(a)} is a tail-dependence matrix.
\item[(ii)] For $\alpha\in[0,1]$, if $X_i=\prod_{j=i}^{i+d-1}
Z_j$, $i=1,\ldots,d$, for independent $\mathrm{B}(1,\alpha)$ random variables
$Z_1,\ldots,Z_{2d-1}$, then $\Lambda$ is an AR(1) matrix with
parameter $\alpha$;
that is, \textup{(b)} is a tail-dependence matrix.
\item[(iii)] For $\alpha\in[0,1/2]$, if $X_i=\id_{\{Z\in
[(i-1)(1-\alpha), (i-1)(1-\alpha)+1]\}}$, $i=1,\ldots,d$, for
$Z\sim\mathrm{U}[0,d]$, then $\Lambda$ is an MA(1) matrix with
parameter $\alpha$;
that is, \textup{(c)} is a tail-dependence matrix.
\end{longlist}
\end{proposition}
\begin{pf}
We have seen in the proof of Theorem~\ref{thm:main:characterization}
that if
$\E[X_1]=\cdots=\E[X_d]=p$, then $\mathbf{Y}$ defined through
\eqref{smodel2} has
tail-dependence matrix $\E[\mathbf{X}\mathbf{X}^\top]/p$. Write
$\Lambda=(\lambda_{ij})_{d\times
d}$ and note that $\lambda_{ii}=1$, $i=1,\ldots,d$, is always guaranteed.
\begin{longlist}[(iii)]
\item[(i)] For $i\neq j$, we have that $\E[X_iX_j]=\alpha^2$ and thus
$\lambda_{ij}=\alpha^2/\alpha=\alpha$. This shows that
$\Lambda$ is an equicorrelation matrix with parameter $\alpha$.
\item[(ii)] For $i<j$, we have that
\begin{eqnarray*}
\E[X_iX_j]&=&\E \Biggl[ \prod
_{k=i}^{i+d-1} Z_k \prod
_{l=j}^{j+d-1} Z_l \Biggr]=\E \Biggl[ \prod
_{k=i}^{j-1} Z_k \Biggr]\E \Biggl[
\prod_{k=j}^{i+d-1} Z_k \Biggr]\E
\Biggl[ \prod_{k=i+d}^{j+d-1} Z_k
\Biggr]
\\
&=&\alpha^{j-i} \alpha^{i+d-j}\alpha^{j-i}=
\alpha^{j-i+d}
\end{eqnarray*}
and $\E[X_i]=\E[X_i^2]=\alpha^d$. Hence,
$\lambda_{ij}=\alpha^{j-i+d}/\alpha^{d}=\alpha^{j-i}$ for $i< j$. By
symmetry, $\lambda_{ij}= \alpha^{ |i-j| }$ for $i\ne j$. Thus,
$\Lambda$
is an AR(1) matrix with parameter $\alpha$.
\item[(iii)] For $i<j$, note that $2(1-\alpha)\ge1$, so
\begin{eqnarray*}
\E[X_iX_j]&=&\p\bigl(Z\in\bigl[(j-1) (1-\alpha), (i-1)
(1-\alpha)+1\bigr]\bigr)
\\
&=&\id_{\{j=i+1\}} \p\bigl(Z\in\bigl[i(1-\alpha),(i-1) (1-\alpha)+1\bigr]
\bigr)=\id_{\{
j=i+1\}}\frac{\alpha} d
\end{eqnarray*}
and $\E[X_i]=\E[X_i^2]=\frac{1}d$. Hence, $\lambda_{ij}=\alpha
\id_{\{j-i=1\}}$ for $i<j$. By symmetry, $\lambda_{ij}=
\alpha\id_{\{|i-j|=1\}}$ for $i\ne j$. Thus, $\Lambda$ is an MA(1)
matrix with
parameter $\alpha$.\quad\qed 
\end{longlist}
\noqed\end{pf}

\subsection{Advanced tail-dependence models}
Theorem~\ref{thm:main:characterization} gives a characterization of
tail-dependence matrices using
Bernoulli-compatible matrices and \eqref{smodel} provides a compatible
model $\mathbf{Y}$
for any tail-dependence matrix $\Lambda(=\E[\mathbf{X}\mathbf
{X}^\top]/p)$.

It is generally not easy to check whether a given matrix is a
Bernoulli-compatible matrix or a tail-dependence matrix; see also
Remark~\ref{rem:Bd}. Therefore, we now study the following question.
%
\begin{question}\label{Q:model:for:Y}
How can we construct a broader class of models with flexible dependence
structures and desired tail-dependence matrices?
\end{question}
To enrich our models, we bring random matrices with Bernoulli entries
into play. For $d, m\in\N$, let
\[
\V_{d\times m}= \Biggl\{X=(X_{ij})_{d\times m}: \p\bigl(X\in
\{0,1\} ^{d\times m}\bigr)=1, \sum_{j=1}^m
X_{ij}\le1, i=1,\ldots,d \Biggr\},
\]
that is, $\V_{d\times m}$ is the set of $d\times m$ random matrices
supported in
$\{0,1\}^{d\times m}$ with each row being \emph{mutually exclusive}; see
Dhaene and Denuit \cite{DD99}. Furthermore,
we introduce a transformation $\L$ on the set of square matrices, such
that, for
any $i,j=1,\ldots,d$, the $(i,j)$th element $\tilde{b}_{ij}$ of $\L
(B)$ is
given by
%
\begin{equation}
\tilde{b}_{ij}=\cases{ b_{ij},&\quad$\mbox{if } i\neq j$,
\vspace*{2pt}
\cr
1,&\quad$\mbox{if } i=j$;}
\end{equation}
that is, $\L$ adjusts the diagonal entries of a matrix to be 1, and
preserves all the other entries. For a set $S$ of square matrices, we set
$\L(S)=\{\L(B):B\in S\}$. We can now address Question~\ref{Q:model:for:Y}.
%
\begin{theorem}[(A class of flexible models)]\label{thm:main:tdm:Y}
Let $\mathbf{U}\sim C^{\mathbf{U}}$ for an $m$-dimensional copula
$C^{\mathbf{U}}$ with tail-dependence
matrix $\Lambda$ and let $\mathbf{V}\sim C^{\mathbf{V}}$ for a
$d$-dimensional copula $C^{\mathbf{V}}$
with tail-dependence matrix $I_d$. Furthermore, let $X \in\V_{d\times
m}$ such that $X,\mathbf{U},\mathbf{V}$ are independent and let
%
\begin{equation}
\mathbf{Y}=X\mathbf{U}+\mathbf{Z}\circ\mathbf{V},\label{eq:stoch:rep}
\end{equation}
where $\mathbf{Z}=(Z_1,\ldots,Z_d)$ with $Z_i=1-\sum_{k=1}^mX_{ik}$,
$i=1,\ldots,d$.
Then $\mathbf{Y}$ has tail-dependence matrix $\Gamma=\L(\E[X\Lambda
X^\top])$.
\end{theorem}
\begin{pf}
Write $X=(X_{ij})_{d\times m}$, $\mathbf{U}=(U_1,\ldots,U_m)$,
$\mathbf{V}=(V_1,\ldots,V_d)$, $\Lambda=(\lambda_{ij})_{d\times d}$ and
$\mathbf{Y}=(Y_1,\ldots,Y_d)$. Then, for all $i=1,\ldots,d$,
\begin{eqnarray*}
Y_i=\sum_{k=1}^m
X_{ik} U_k+Z_iV_i=\cases{
V_i,&\quad$\mbox{if } X_{ik}=0 \mbox{ for all } k=1,\ldots,m,
\mbox{ so } Z_i=1$,\vspace*{2pt}
\cr
U_k,&\quad$\mbox{if }
X_{ik}=1 \mbox{ for some } k=1,\ldots,m, \mbox{ so } Z_i=0$.}
\end{eqnarray*}
Clearly, $\mathbf{Y}$ has $\mathrm{U}[0,1]$ margins. We now calculate the
tail-dependence matrix $\Gamma=(\gamma_{ij})_{d\times d}$ of $Y$ for
$i\neq j$. By our
independence assumptions, we can derive the following results:
\begin{longlist}[(iii)]
\item[(i)]$\p(Y_i\le u, Y_j\le u, Z_i=1,
Z_j=1)=\p(V_i\le u, V_j\le u, Z_i=1, Z_j=1)=C_{ij}^{\mathbf
{V}}(u,u)\p(Z_i=1,Z_j=1)\le
C_{ij}^{\mathbf{V}}(u,u)$, where $C_{ij}^{\mathbf{V}}$ denotes the
$(i,j)$th margin of $C^{\mathbf{V}}$. As
$\mathbf{V}$ has tail-dependence matrix $I_d$, we obtain that
\[
\lim_{u\downarrow0}\frac{1}u\p(Y_i\le u,
Y_j\le u, Z_i=1, Z_j=1)=0.
\]
\item[(ii)]$\p(Y_i\le u, Y_j\le u, Z_i=0, Z_j=1)=\sum_{k=1}^m\p(U_k\le u,
V_j\le u,
X_{ik}=1, Z_j=1)=\sum_{k=1}^m\p(U_k\le u)\p(V_j\le u)\p
(X_{ik}=1,Z_j=1)\le
u^2$, and thus
\[
\lim_{u\downarrow0}\frac{1}u\p(Y_i\le u,
Y_j\le u, Z_i=0, Z_j=1)=0.
\]
Similarly, we obtain that
\[
\lim_{u\downarrow0}\frac{1}u\p(Y_i\le u,
Y_j\le u, Z_i=1, Z_j=0)=0.
\]
\item[(iii)]$\p(Y_i\le u, Y_j\le u, Z_i=0, Z_j=0)=\sum_{k=1}^m\sum_{l=1}^m\p(U_k\le
u, U_l\le u,  X_{ik}=1,
X_{jl}=1)=\sum_{k=1}^m\sum_{l=1}^mC_{kl}^{\mathbf{U}}(u,u)\p(X_{ik}=1,
X_{jl}=1)= \sum_{k=1}^m\sum_{l=1}^mC_{kl}^{\mathbf{U}}(u,u)\*\E
[X_{ik}X_{jl}]$ so that
\begin{eqnarray*}
\lim_{u\downarrow0}\frac{1}u\p(Y_i\le u,
Y_j\le u, Z_i=0, Z_j=0)&=&\sum
_{k=1}^m\sum_{l=1}^m
\lambda_{kl}\E[X_{ik}X_{jl}]\\
&=&\E \Biggl[ \sum
_{k=1}^m\sum_{l=1}^mX_{ik}
\lambda_{kl}X_{jl} \Biggr]
\\
&= &\bigl(\E\bigl[X\Lambda X^\top\bigr] \bigr)_{ij}.
\end{eqnarray*}
\end{longlist}
By the law of total probability, we thus obtain that
\begin{eqnarray*}
\gamma_{ij}&=&\lim_{u\downarrow0}\frac{\p(Y_i\le u, Y_j\le
u)}{u}=\lim
_{u\downarrow0}\frac{\p(Y_i\le u, Y_j\le u, Z_i=0,
Z_j=0)}{u}
\\
&=& \bigl(\E\bigl[X\Lambda X^\top\bigr] \bigr)_{ij}.
\end{eqnarray*}
This shows that $\E[X\Lambda X^\top]$ and $\Gamma$ agree on the off-diagonal
entries. Since $\Gamma\in\T_d$ implies that $\diag(\Gamma)=I_d$,
we conclude
that $\L(\E[X\Lambda X^\top])=\Gamma$.
\end{pf}
%

A special case of Theorem~\ref{thm:main:tdm:Y} reveals an essential difference
between the transition rules of a tail-dependence matrix and a
covariance matrix.
Suppose that for $X\in\V_{d\times m}$, $\E[X]$ is a stochastic
matrix (each row
sums to 1), and $\mathbf{U}\sim C^{\mathbf{U}}$ for an
$m$-dimensional copula $C^{\mathbf{U}}$ with
tail-dependence matrix $\Lambda=(\lambda_{ij})_{d\times d}$. Now we
have that
$Z_i=0$, $i=1,\ldots,d$ in \eqref{eq:stoch:rep}. By Theorem~\ref{thm:main:tdm:Y}, the tail dependence matrix of $\mathbf
{Y}=X\mathbf{U}$ is given by
$\L(\E[X\Lambda X^\top])$. One can check the diagonal terms of the matrix
$\Lambda^*=(\lambda^*_{ij})_{d\times d}=X\Lambda X^\top$ by
\[
\lambda^*_{ii}=\sum_{j=1}^m
\sum_{k=1}^m X_{ik}\lambda
_{kj}X_{ij}=\sum_{k=1}^m
X_{ik}\lambda_{kk}=1,\qquad i=1,\ldots,m.
\]
Hence, the tail-dependence matrix of $\mathbf{Y}$ is indeed $\E
[X\Lambda X^\top]$.

\begin{remark} \label{transition} In summary:
\begin{longlist}[(ii)]
\item[(i)] If an $m$-vector $\mathbf{U}$ has covariance matrix $\Sigma$,
then $X\mathbf{U}$ has
covariance matrix $\E[X \Sigma X^\top]$ for any $d\times m$ random matrix
$X$ independent of $\mathbf{U}$.
\item[(ii)] If an $m$-vector $\mathbf{U}$ has uniform $[0,1]$ margins and
tail-dependence
matrix $\Lambda$, then $X\mathbf{U}$ has tail-dependence matrix $\E
[X\Lambda X^\top]$
for any $X\in\V_{d\times m}$ independent of $\mathbf{U}$ such that
each row of $X$
sums to 1.
\end{longlist}
It is noted that the transition property of tail-dependence matrices is more
restricted than that of covariance matrices.
\end{remark}

The following two propositions consider selected special cases of this
construction which are more straightforward to apply.
%
\begin{proposition}\label{coro-b}
For any $B \in\B_d$ and any $\Lambda\in\T_d$ we have that $\L
(B\circ
\Lambda)\in\T_d$. In particular, $\L(B) \in\T_d$, and hence
$\L(\B_d)\subseteq\T_d$.
\end{proposition}
\begin{pf}
Write $B=(b_{ij})_{d\times d}=\E[\mathbf{W}\mathbf{W}^\top]$ for
some $\mathbf{W}=(W_1,\ldots,W_d)\in
\V_d$ and consider $X=\diag(\mathbf{W})\in\V_{d\times d}$. As in
the proof of
Theorem~\ref{thm:main:tdm:Y} (and with the same notation), it follows
that for $i\ne
j$, $\gamma_{ij}=\E[X_{ii}\lambda_{ij}X_{jj}]=\E[W_iW_j\lambda
_{ij}]$. This
shows that $\E[X\Lambda X^\top]=\E[\mathbf{W}\mathbf{W}^\top\circ
\Lambda]$ and $B\circ\Lambda$
agree on off-diagonal entries. Thus, $\L(B\circ\Lambda)=\Gamma\in
\T_d$. By
taking $\Lambda=(1)_{d\times d}$, we obtain $\L(B)\in T_d$.
\end{pf}

The following proposition states a relationship between substochastic
matrices and
tail-dependence matrices. To this end, let
\[
\Q_{d\times m}= \Biggl\{Q=(q_{ij})_{d\times m}: \sum
_{j=1}^m q_{ij}\le 1, q_{ij}
\ge0, i=1,\ldots,d, j=1,\ldots,m \Biggr\},
\]
that is, $\Q_{d\times m}$ is the set of $d\times m$ \emph{(row) substochastic
matrices}; note that the expectation of a random matrix in $\mathcal
V_{d\times m}$
is a substochastic matrix.
%
\begin{proposition}
For any $Q\in\Q_{d\times m}$ and any $\Lambda\in\T_m$, we have
that $\L(Q
\Lambda Q^\top)\in\T_d$. In particular, $\L(QQ^\top) \in\T_d$
for all $Q\in
\Q_{d\times m}$ and $\L(\bd p \bd p^\top) \in\T_d$ for all $\bd
p\in[0,1]^d$.
\end{proposition}
\begin{pf}
Write $Q=(q_{ij})_{d\times m}$ and let $X_{ik}=\id_{\{Z_{i}\in
[ \sum_{j=1}^{k-1}q_{ij}, \sum_{j=1}^{k}q_{ij})\}}$ for independent
$Z_i\sim\mathrm{U}[0,1]$, $i=1,\ldots,d$, $k=1,\ldots,m$. It is
straightforward
to see that $\E[X]=Q$, $X\in\V_{d\times m}$ with independent rows, and
$\sum_{k=1}^m X_{ik}\le1$ for $i=1,\ldots,d$, so $X\in\V_{d\times
m}$. As in
the proof of Theorem~\ref{thm:main:tdm:Y} (and with the same
notation), it follows that
for $i\ne j$,
\[
\gamma_{ij}= \sum_{l=1}^m\sum
_{k=1}^m \E[X_{ik}]
\E[X_{jl}]\lambda_{kl}= \sum_{l=1}^m
\sum_{k=1}^m q_{ik}q_{jl}
\lambda_{kl}.
\]
This shows that $Q\Lambda Q^\top$ and $\Gamma$ agree on off-diagonal
entries, so $\L(Q\Lambda Q^\top)=\Gamma\in\T_d$. By taking
$\Lambda=I_d$, we obtain $\L(QQ^\top)\in T_d$. By taking $m=1$, we obtain
$\L(\bd p \bd p^\top) \in\T_d$.
\end{pf}

\subsection{Corresponding copula models}
In this section, we derive the copulas of \eqref{smodel} and~\eqref{eq:stoch:rep} which are able to produce tail-dependence matrices
$\E[\mathbf{X}\mathbf{X}^\top]/p$ and $\L(\E[X\Lambda X^\top])$
as stated in
Theorems~\ref{thm:main:characterization} and \ref{thm:main:tdm:Y},
respectively. We first address the former.
%
\begin{proposition}[{[Copula of \eqref{smodel}]}]
Let $\mathbf{X}\in\V_d$, $\E[\mathbf{X}]=(p,\ldots,p)\in(0,1]^d$.
Furthermore, let
$U,V\sim\mathrm{U}[0,1]$, $U,V,\mathbf{X}$ be independent and
\[
\mathbf{Y}=\mathbf{X}pU+(\mathbf{1}-\mathbf{X}) \bigl(p+(1-p)V\bigr).
\]
Then the copula $C$ of $\mathbf{Y}$ at $\mathbf{u}=(u_1,\ldots,u_d)$
is given by
\[
C(\mathbf{u})=\sum_{\mathbf{i}\in\{0,1\}^d}\min \biggl\{
\frac{\min_{r:i_r=1}\{u_r\}}{p},1 \biggr\} \max \biggl\{\frac{\min_{r:i_r=0}\{
u_r\}-p}{1-p},0 \biggr\} \p(
\mathbf{X}=\mathbf{i}),
\]
with the convention $\min\varnothing=1$.
\end{proposition}
\begin{pf}
By the law of total probability and our independence assumptions,
\begin{eqnarray*}
C(\mathbf{u})&=& \sum_{\mathbf{i}\in\{0,1\}^d}\p(\mathbf{Y}\le
\mathbf{u}, \mathbf{X}=\mathbf{i})
\\[-2pt]
&=& \sum_{\mathbf{i}\in\{0,1\}^d}\p\Bigl(pU\le\min_{r:i_r=1}
\{u_r\}, p+(1-p)V\le\min_{r:i_r=0}\{u_r\},
\mathbf{X}=\mathbf{i}\Bigr)
\\[-2pt]
&= &\sum_{\mathbf{i}\in\{0,1\}^d}\p \biggl(U\le\frac{\min_{r:i_r=1}\{u_r\}}{p}
\biggr) \p \biggl(V\le\frac{\min_{r:i_r=0}\{
u_r\}-p}{1-p} \biggr) \p(\mathbf{X}=\mathbf{i});
\end{eqnarray*}
the claim follows from the fact that $U,V\sim\mathrm{U}[0,1]$.
\end{pf}

For deriving the copula of \eqref{eq:stoch:rep}, we need to introduce some
notation; see also Example~\ref{ex:special:cases} below. In the following
theorem, let $\mathrm{supp}(X)$ denote the support of~$X$. For a vector
 $\mathbf{u}=(u_1,\ldots,u_d)\in[0,1]^d$ and a matrix
$A=(A_{ij})_{d\times m}\in\mathrm{supp}(X)$, denote by $A_i$ the sum
of the
$i$th row of $A$, $i=1,\ldots,d$, and let
$\mathbf{u}_A=(u_1\id_{\{A_1=0\}}+\id_{\{A_1=1\}},\ldots,u_d\id_{\{
A_d=0\}}+\id_{\{A_d=1\}})$,
and $\mathbf{u}_A^*=(\min_{r:A_{r1}=1}\{u_r\},\ldots,
\min_{r:A_{rm}=1}\{u_r\})$,
where $\min\varnothing=1$.

%
\begin{proposition}[{[Copula of \eqref{eq:stoch:rep}]}]
Suppose that the setup of Theorem~\ref{thm:main:tdm:Y} holds. Then the copula
$C$ of $\mathbf{Y}$ in \eqref{eq:stoch:rep} is given by
%
\begin{equation}
C(\mathbf{u})=\sum_{A\in\mathrm{supp}(X)} C^V(
\mathbf{u}_A) C^U\bigl(\mathbf{u}_A^*\bigr)
\p(X=A).\label{eq:cop:Y}
\end{equation}
%
\end{proposition}
\begin{pf}
By the law of total probability, it suffices to verify that $\p
(\mathbf{Y}\le
\mathbf{u} | X=A)= C^{\mathbf{V}}(\mathbf{u}_A)C^{\mathbf
{U}}(\mathbf{u}_A^*)$. This can be seen from
\begin{eqnarray*}
&&\p(\mathbf{Y}\le\mathbf{u} | X=A)
\\[-2pt]
&&\qquad=\p \Biggl( \sum_{k=1}^m
A_{jk}U_k+(1-A_j)V_j\le
u_j, j=1,\ldots ,d \Biggr)
\\[-2pt]
&&\qquad=\p(U_k \id_{\{A_{jk}=1\}} \le u_j, V_j
\id_{\{A_j=0\}}\le u_j, j=1,\ldots,d,~k=1,\ldots,m)
\\[-2pt]
&&\qquad=\p \Bigl(U_k \le\min_{r:A_{rk}=1}\{u_r\},
V_j\le u_j \id_{\{A_j=0\}}+\id_{\{A_j=1\}},\\[-2pt]
&&\qquad\quad j=1,\ldots,d, k=1,\ldots,m \Bigr)
\\[-2pt]
&&\qquad=\p \Bigl(U_k \le\min_{r:A_{rk}=1}\{u_r
\}, k=1,\ldots,m \Bigr) \p (V_j\le u_j
\id_{\{A_j=0\}}+\id_{\{A_j=1\}},\\[-2pt]
&&\qquad\quad j=1,\ldots,d )
\\[-2pt]
&&\qquad=C^{\mathbf{U}}\bigl(\mathbf{u}_A^*\bigr) C^{\mathbf{V}}({
\mathbf{u}}_A).%
\end{eqnarray*}
\upqed\end{pf}

As long as $C^{\mathbf{V}}$ has tail-dependence matrix $I_d$, the
tail-dependence matrix of
$\mathbf{Y}$ is not affected by the choice of $C^{\mathbf{V}}$. This
theoretically provides more
flexibility in choosing the body of the distribution of $\mathbf{Y}$
while attaining a
specific tail-dependence matrix. Note, however, that this also depends
on the
choice of $X$; see the following example where we address special cases which
allow for more insight into the rather abstract construction~\eqref{eq:cop:Y}.
%
\begin{example}\label{ex:special:cases}
%
1. For $m=1$, the copula $C$ in \eqref
{eq:cop:Y} is given by
%
\begin{equation}
C(\mathbf{u})=\sum_{\mathbf{A}\in\{0,1\}^d} C^{\mathbf{V}}(\mathbf
{u}_{\mathbf{A}})C^{\mathbf{U}}\bigl(\mathbf{u}_{\mathbf{A}}^*\bigr)\p (
\mathbf{X}=\mathbf{A});\label{eq:case:m1}
\end{equation}
note that $X,A$ in equation~(\ref{eq:cop:Y}) are indeed vectors in this
case. For $d=2$, we obtain
\begin{eqnarray*}
C(u_1,u_2)&=&M(u_1,u_2)\p\biggl(
\mathbf{X}=\pmatrix{1\cr 1}\biggr)+C^{\mathbf
{V}}(u_1,u_2)\p\biggl(
\mathbf{X}=\pmatrix{0\cr 0}\biggr)
\\
&&{}+\Pi(u_1,u_2)\p\biggl(\mathbf{X}=\pmatrix{1\cr 0}
\mbox{ or }\mathbf{X}=\pmatrix{0\cr 1}\biggr),
\end{eqnarray*}
and therefore a mixture of the Fr\'echet--Hoeffding upper bound
$M(u_1,u_2)$ $=\min\{u_1,u_2\}$, the copula $C^{\mathbf{V}}$ and the
independence copula
$\Pi(u_1,u_2)=u_1u_2$. If $\p\bigl(\mathbf{X}=\bigl({0\atop 0}\bigr)\bigr)=0$ then $C$
is simply a mixture of $M$ and $\Pi$ and does not depend on $\mathbf
{V}$ anymore.

Now consider the special case of \eqref{eq:case:m1} where $\mathbf
{V}$ follows the
$d$-dimensional independence copula $\Pi(\mathbf{u})=\prod_{i=1}^d u_i$
and $\mathbf{X}=(X_1,\ldots,X_{d-1},1)$ is such that at most one of
$X_1,\ldots,X_{d-1}$ is 1 [each randomly with probability $0\le\alpha
\le
1/(d-1)$ and all are simultaneously 0 with probability
$1-(d-1)\alpha$]. Then, for all $\mathbf{u}\in[0,1]^d$, $C$ is given by
%
\begin{equation}
C(\mathbf{u}) 
=\alpha\sum
_{i=1}^{d-1} \Biggl(\min\{u_i,u_d
\}\prod_{j=1, j\ne
i}^{d-1}u_j \Biggr) +
\bigl(1-(d-1)\alpha\bigr) \prod_{j=1}^{d}u_j.\label
{ex:special:cases:FD}
\end{equation}
%
This copula is a conditionally independent multivariate Fr\'echet copula
studied in Yang et al. \cite{YQW09}. This example will be revisited in
Section~\ref{sec:ex:qrm}; see also the left-hand side of
Figure~\ref{fig:Gamma} below.

\begin{figure}[b]

\includegraphics{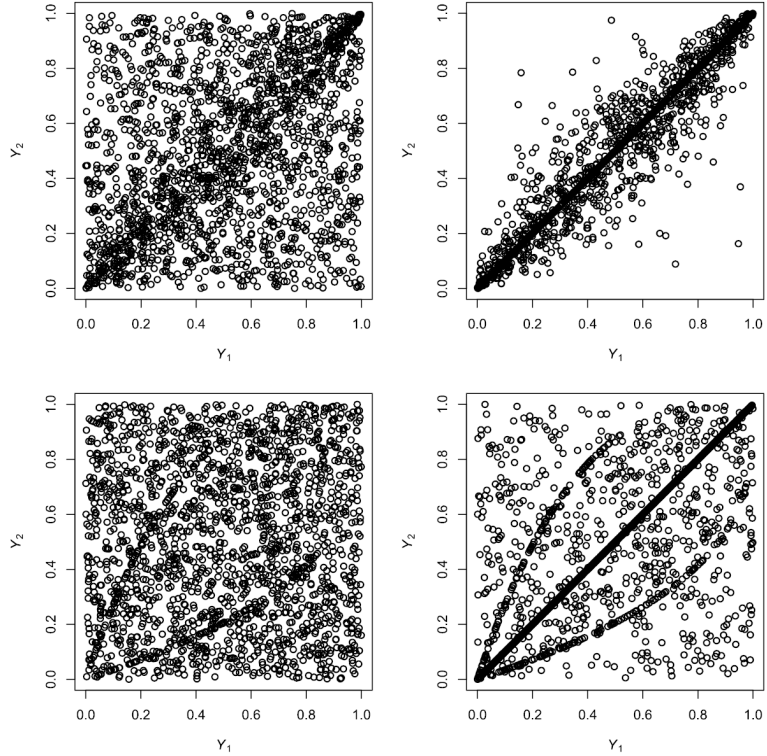}

\caption{Scatter plots of 2000 samples from $\mathbf{Y}$ for $\mathbf
{V}\sim\Pi$ and $\mathbf{U}$
following a bivariate ($m=2$) $t_3$ copula with Kendall's tau equal to 0.75
(top row) or a survival Marshall--Olkin copula with parameters
$\alpha_1=0.25,\alpha_2=0.75$ (bottom row). For the plots on the left-hand
side, the number of rows of $X$ with one 1 are randomly chosen among
$\{0,1,2\ (=d)\}$, the corresponding rows and columns are then randomly
selected among $\{1,2\ (=d)\}$ and $\{1,2\ (=m)\}$, respectively. For the plots
on the right-hand side, $X$ is drawn from a multinomial distribution with
probabilities 0.5 and 0.5 such that each row contains precisely one~1.}
\label{fig:Y:m2}
\end{figure}

2. For $m=2$, $d=2$, we obtain
%
\begin{eqnarray}
\label{ex:C:d2:m2} C(u_1,u_2)&=&M(u_1,u_2)
\p\biggl(X=\mat{1} {0} {1} {0} \mbox{ or } X=\mat {0} {1} {0} {1}\biggr)
\nonumber
\\
&&{}+C^{\mathbf{U}}(u_1,u_2)\p\biggl(X=\mat {1}
{0} {0} {1}\biggr)+C^{\mathbf{U}}(u_2,u_1)\p\biggl(X=\mat{0} {1}
{1} {0}\biggr)
\nonumber
\\[-8pt]
\\[-8pt]
\nonumber
&&{}+C^{\mathbf{V}}(u_1,u_2)\p\biggl(X=\mat{0}
{0} {0} {0}\biggr)\\
&&{}+\Pi (u_1,u_2)\p\biggl(X=\mat{0} {0} {1} {0}
\mbox{ or } \mat{0} {0} {0} {1} \mbox { or } \mat{1} {0} {0} {0}
\mbox{ or } \mat{0} {1}
{0} {0}\biggr).\hspace*{-25pt}\nonumber
\end{eqnarray}
Figure~\ref{fig:Y:m2} shows samples of size 2000 from \eqref
{ex:C:d2:m2} for
$\mathbf{V}\sim\Pi$ and two different choices of $\mathbf{U}$ (in
different rows) and $X$
(in different columns). From Theorem~\ref{thm:main:tdm:Y}, we obtain
that the
off-diagonal entry $\gamma_{12}$ of the tail-dependence matrix $\Gamma
$ of
$\mathbf{Y}$ is given by
\[
\gamma_{12}=p_{(1,2)(1,1)}+p_{(1,2)(2,2)}+\lambda
_{12}(p_{(1,2)(2,1)}+p_{(1,2)(1,2)}),
\]
where $\lambda_{12}$ is the off-diagonal entry of the tail-dependence
matrix $\Lambda$ of $\mathbf{U}$.
\end{example}

\subsection{An example from risk management practice}\label{sec:ex:qrm}
Let us now come back to problem~\eqref{eq:Gammad} which motivated our research
on tail-dependence matrices. From a practical point of view, the
question is
whether it is possible to find one financial position, which has tail-dependence
coefficient $\alpha$ with each of $d-1$ tail-independent financial risks
(assets). Such a construction can be interesting for risk management purposes,
for example, in the context of hedging.\looseness=-1

Recall problem~\eqref{eq:Gammad}:

\begin{quote}
\textit{For which $\alpha\in[0,1]$ is the matrix
%
\begin{eqnarray}
\Gamma_d(\alpha)=\pmatrix{ 1 & 0 & \cdots& 0 & \alpha\vspace*{2pt}
\cr
0 & 1 & \cdots& 0 & \alpha\vspace*{2pt}
\cr
\vdots& \vdots& \ddots& \vdots&
\vdots\vspace*{2pt}
\cr
0 & 0 & \cdots& 1 & \alpha\vspace*{2pt}
\cr
\alpha& \alpha&
\cdots& \alpha& 1 } \label{eq:Gammadrec}
\end{eqnarray}
a matrix of pairwise (either lower or upper) tail-dependence coefficients?}
\end{quote}

Based on the Fr\'echet--Hoeffding bounds, it follows from
Joe \cite{J97}, Theorem~3.14, that for $d=3$ (and thus also $d>3$),
$\alpha$ has to
be in $[0,1/2]$; however, this is not a sufficient condition for
$\Gamma_d(\alpha)$ to be a tail-dependence matrix. The following
proposition not
only gives an answer to \eqref{eq:Gammadrec} by providing necessary and
sufficient such conditions, but also provides, by its proof, a
compatible model
for $\Gamma_d(\alpha)$.
%
\begin{proposition}\label{prop:example}
$\Gamma_d(\alpha) \in\T_d$ if and only if $0\le\alpha\le1/(d-1)$.
\end{proposition}
\begin{pf}
The if-part directly follows from Corollary~\ref{coro:ddm}. We provide
a constructive proof based on Theorem~\ref{thm:main:tdm:Y}.
Suppose that $0\le\alpha\le1/(d-1)$. Take a partition
$\{\Omega_1,\ldots,\Omega_{d}\}$ of the sample space $\Omega$ with
$\p(\Omega_i)=\alpha$, $i=1,\ldots,d-1$, and let
$\mathbf{X}=(\id_{\Omega_1},\ldots,\id_{\Omega_{d-1}}, 1)\in\V
_d$. It is
straightforward to see that\looseness=-1
\[
\E\bigl[\mathbf{X}\mathbf{X}^\top\bigr]=\pmatrix{ \alpha& 0 & \cdots&
0 & \alpha\vspace*{2pt}
\cr
0 & \alpha& \cdots& 0 & \alpha\vspace*{2pt}
\cr
\vdots&
\vdots& \ddots& \vdots&\vdots \vspace*{2pt}
\cr
0 & 0 & \cdots& \alpha& \alpha
\vspace*{2pt}
\cr
\alpha& \alpha& \cdots& \alpha& 1 }.
\]
By Proposition~\ref{coro-b}, $\Gamma_d(\alpha)=\L(\E[\mathbf
{X}\mathbf{X}^\top])\in\T_d$.

For the only if part, suppose that $\Gamma_d(\alpha) \in\T_d$; thus
$\alpha\ge0$. By Theorem~\ref{thm:main:characterization}, $\Gamma
_d(\alpha)\in
\B^I_d$. By the definition of $\B^I_d$, $\Gamma_d(\alpha)=B_d/p$
for some
$p\in(0,1]$ and 
a Bernoulli-compatible matrix $B_d$. 
Therefore,
\[
p\Gamma_d(\alpha)=\pmatrix{ p & 0 & \cdots& 0 & p\alpha\vspace*{2pt}
\cr
0 & p & \cdots& 0 & p\alpha\vspace*{2pt}
\cr
\vdots& \vdots& \ddots& \vdots&
\vdots \vspace*{2pt}
\cr
0 & 0 & \cdots& p & p\alpha\vspace*{2pt}
\cr
p\alpha& p
\alpha& \cdots& p\alpha& p }
\]
is a compatible Bernoulli matrix, so $p\Gamma_d(\alpha)\in\B_d$. Write
$p\Gamma_d(\alpha)=\E[\mathbf{X}\mathbf{X}^\top]$ for some
$\mathbf{X}=(X_1,\ldots,X_d)\in\V_d$. It
follows that $\p(X_i=1)=p$ for $i=1,\ldots, d$, $\p(X_iX_j=1)=0$ for
$i\ne j$,
$i,j=1,\ldots,d-1$ and $\p(X_iX_d=1)=p\alpha$ for $i=1,\ldots,d-1$.
Note that
$\{X_iX_d=1\}$, $i=1,\ldots,d-1$, are almost surely disjoint since $\p
(X_iX_j=1)=0$ for
$i\ne j$, $i,j=1,\ldots,d-1$. As a consequence,
\[
p=\p(X_d=1)\ge\p \Biggl( \bigcup_{i=1}^{d-1}
\{X_iX_d=1\} \Biggr)=\sum_{i=1}^{d-1}
\p(X_iX_d=1)=(d-1)p\alpha,
\]
and thus $(d-1)\alpha\le1$.
\end{pf}

\begin{figure}

\includegraphics{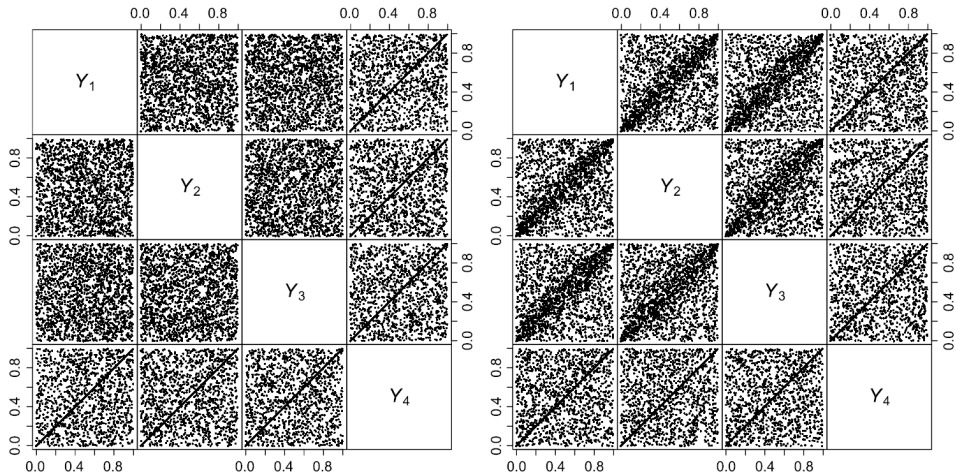}

\caption{Pairs plot of 2000 samples from $\mathbf{Y}\sim C$ which
produces the tail
dependence matrix $\Gamma_4(1/3)$ as given by \protect\eqref
{eq:Gammad}. On the
left-hand side, $\mathbf{V}\sim\Pi$ [$\alpha$ determines how much
weight is on the
diagonal for pairs with one component being $Y_4$; see
\protect\eqref{ex:special:cases:FD}] and on the right-hand side,
$\mathbf{V}$ follows a Gauss
copula with parameter chosen such that Kendall's tau equals 0.8.}
\label{fig:Gamma}
\end{figure}

It follows from the proof of Theorem~\ref{thm:main:tdm:Y} that for
$\alpha\in[0,1/(d-1)]$, a compatible copula model with
tail-dependence matrix
$\Gamma_d(\alpha)$ can be constructed as follows. Consider a partition
$\{\Omega_1,\ldots,\Omega_{d}\}$ of the sample space $\Omega$ with
$\p(\Omega_i)=\alpha$, $i=1,\ldots,d-1$, and let
$\mathbf{X}=(X_1,\ldots,X_d)=(\id_{\Omega_1},\ldots,\id_{\Omega
_{d-1}}, 1)\in\V_d$; note
that $m=1$ here. Furthermore, let $\mathbf{V}$ be as in Theorem~\ref
{thm:main:tdm:Y}, $U\sim\mathrm{U}[0,1]$
and $U,\mathbf{V},\mathbf{X}$ be independent. Then
\[
\mathbf{Y}=\bigl(UX_1+(1-X_1)V_1,
\ldots,UX_{d-1}+(1-X_{d-1})V_{d-1},U\bigr)
\]
has tail-dependence matrix $\Gamma_d(\alpha)$. Example~\ref
{ex:special:cases},
part~1 provides the copula $C$ of $\mathbf{Y}$
in this case. It is
also straightforward to verify from this copula that $\mathbf{Y}$ has
tail-dependence matrix
$\Gamma_d(\alpha)$. 
Figure~\ref{fig:Gamma} displays pairs plots of 2000 realizations of
$\mathbf{Y}$ for
$\alpha=1/3$ and two different copulas for $\mathbf{V}$.

\begin{remark}
Note that $\Gamma_d(\alpha)$ is not positive semidefinite if and only if
$\alpha>1/\sqrt{d-1}$. For $d<5$, element-wise nonnegative and positive
semidefinite matrices are completely positive; see Berman and
Shaked-Monderer \cite{BS03}, Theorem~2.4. Therefore, $\Gamma_3(2/3)$ is completely positive. However, it is
not in
$\T_3$. It indeed shows that the class of completely positive matrices with
diagonal entries being 1 is strictly larger than $\T_d$.
\end{remark}

\section{Conclusion and discussion}\label{sec:con}
Inspired by the question whether a given matrix in $[0,1]^{d\times d}$
is the
matrix of pairwise tail-dependence coefficients of a $d$-dimensional random
vector, we introduced the tail-dependence compatibility problem. It
turns out that this
problem is closely related to the Bernoulli-compatibility problem which
we also
addressed in this paper and which asks when a given matrix in
$[0,1]^{d\times d}$
is a Bernoulli-compatible matrix (see Question~\ref{Q:B:compatibility} and
Theorem~\ref{thm:bern}). As a main finding, we characterized tail-dependence
matrices as precisely those square matrices with diagonal entries being
1 which
are Bernoulli-compatible matrices multiplied by a constant (see
Question~\ref{Q:characterization} and
Theorem~\ref{thm:main:characterization}). Furthermore, we presented
and studied
new models (see, e.g., Question~\ref{Q:model:for:Y} and
Theorem~\ref{thm:main:tdm:Y}) which provide answers to several
questions related
to the tail-dependence compatibility problem.

The study of compatibility of tail-dependence matrices is mathematically
different from that of covariances matrices. Through many technical
arguments in
this paper, the reader may have already realized that the
tail-dependence matrix
lacks a linear structure which is essential to covariance matrices
based on
tools from linear algebra. For instance, let $\mathbf{X}$ be a
$d$-random vector with
covariance matrix $\Sigma$ and tail-dependence matrix $\Lambda$, and
$A$ be an
$m\times d$ matrix. The covariance matrix of $A\mathbf{X}$ is simply
given by $A\Sigma
A^\top$; however, the tail-dependence matrix of $A\mathbf{X}$ is
generally not explicit
(see Remark~\ref{transition} for special cases). This lack of
linearity can also
help to understand why tail-dependence matrices are realized by models
based on Bernoulli
vectors as we have seen in this paper, in contrast to covariance
matrices which
are naturally realized by Gaussian (or generally, elliptical) random
vectors. The latter
have a linear structure, whereas Bernoulli vectors do not. It is not
surprising that most classical techniques in linear algebra such as matrix
decomposition, diagonalization, ranks, inverses and determinants are
not very
helpful for studying the compatibility problems we address in this paper.

Concerning future research, an interesting open question is how one can
(theoretically or numerically) determine whether a given arbitrary nonnegative,
square matrix is a tail-dependence or Bernoulli-compatible matrix. To
the best
of our knowledge there are no corresponding algorithms available.
Another open
question concerns the compatibility of other matrices of pairwise
measures of
association such as rank-correlation measures (e.g., Spearman's rho or Kendall's
tau); see \cite{EMS02}, Section~6.2. Recently, \cite{FSS14} and \cite
{SBS15} studied the
concept of tail-dependence functions of stochastic processes. Similar results
to some of our findings were found in the context of max-stable processes.

From a practitioner's point-of-view, it is important to point out
limitations of
using tail-depen\-dence matrices in quantitative risk management and other
applications. One possible such limitation is the statistical
estimation of
tail-dependence matrices since, as limits, estimating tail dependence
coefficients from data is nontrivial (and typically more complicated than
estimation in the body of a bivariate distribution).

After presenting the results of our paper at the conferences ``Recent
Developments in Dependence Modelling with Applications in Finance and Insurance---2nd Edition, Brussels, May 29, 2015''
and ``The 9th International Conference
on Extreme Value Analysis, Ann Arbor, June 15--19, 2015,'' the references
\cite{FSS14} and \cite{SBS15} were brought to our attention (see also
Acknowledgments below). In these papers, a very related problem is
treated, be
it from a different, more theoretical angle, mainly based on the theory of
max-stable and Tawn--Molchanov processes as well as results for
convex-polytopes. For instance, our Theorem~\ref{thm:main:characterization} is similar to Theorem~6(c) in
\cite{FSS14}.

\section*{Acknowledgments}
We would like to thank Federico Degen for raising this interesting question
concerning tail-dependence matrices and Janusz Milek for relevant discussions.
We would also like to thank the Editor, two anonymous referees,
Tiantian Mao,
Don McLeish, Johan Segers, Kirstin Strokorb and Yuwei Zhao for valuable
input on an early
version of the paper.

%

\printaddresses
\end{document}